\documentclass[11pt,a4paper]{article}
\usepackage[left=26mm,right=26mm,top=23mm,bottom=23mm]{geometry}
\usepackage[T1]{fontenc}
\usepackage{lmodern,microtype,amsmath,amssymb,amsthm,mathtools}
\usepackage[hidelinks]{hyperref}
\hypersetup{pdftitle={A proof of the toughness conjecture for Kneser graphs}}
\newtheorem{theorem}{Theorem}[section]
\newtheorem{lemma}[theorem]{Lemma}

\theoremstyle{definition}
\newcommand{\KG}{\mathrm{KG}}

\title{A proof of the toughness conjecture for Kneser graphs}

\author{Guang Li, Wenqian Zhang\\
{\small School of Mathematics and Statistics, Shandong University of Technology}\\
{\small Zibo, Shandong 255000, P.R. China}
}
\date{}
\begin{document}
\maketitle
\footnotetext{Corresponding author: Wenqian Zhang (zhangwq@pku.edu.cn)}
\begin{abstract}
As a measure of the connectivity of graphs, Chv\'atal \cite{Chvatal} introduced the
concept of toughness. For integers $k\geq2$ and $n\geq2k+1$, the Kneser graph $\KG(n,k)$ has all the $k$-subsets of $[n]=\{1,\ldots,n\}$ as its vertices, and two vertices are adjacent if they are disjoint. Cioab\u{a} and Wong~\cite{CW} determined the toughness of $\KG(n,2)$. Park et al.~\cite{Park} determined the toughness of $\KG(n,k)$ when $k=3,4$ and $k\geq5$ for large $n$. They \cite{Park} proposed a conjecture for $k\geq5$ and other values of $n$. In this paper, we resolve this conjecture.
\end{abstract}

\section{Introduction}

All graphs considered in this paper are finite, undirected and simple.  For a graph $G$, the vertex set and edge set of $G$ are denoted by $V(G)$ and $E(G)$, respectively.  Let $\alpha(G)$ denote the independence number of $G$. Let $c(G)$ denote the number of components of $G$. For connected $G$, a \emph{cut} is any $S\subseteq V(G)$ such that $c(G-S)\ge2$, where $G-S$ denotes the subgraph formed by deleting the vertices in $S$. Observe that deleting $s\geq1$ vertices from a Hamilton cycle leaves at most $s$
components. To measure the connectivity of graphs, Chv\'atal \cite{Chvatal} introduced the
notion of toughness. For a connected
noncomplete graph $G$, the \emph{toughness} is defined as 
\[
       t(G)=\min_{S:\,c(G-S)\ge2}\frac{|S|}{c(G-S)}.
\]
A cut $S$ attaining the minimum is called a \emph{tough set}. It is usually difficult to determine the exact value of toughness of graphs (see \cite{BauerHakimiSchmeichel1990}).

 Let $G$ be a graph with vertices $v_{1},v_{2},\ldots,v_{N}$. The {\em adjacency matrix}  of $G$ is defined as $A(G)=(a_{ij})$, where $a_{ij}=1$ if $v_{i}$ is adjacent to $v_{j}$, and $a_{ij}=0$ otherwise. The eigenvalues of $A(G)$ are denoted by $\lambda_{1}(G)\geq \lambda_{2}(G)\geq\cdots\geq \lambda_{N}(G)$. The {\em Laplacian matrix} of $G$ is defined as $L(G)=D(G)-A(G)$, where $D(G)$ is the diagonal matrix of vertex-degrees of $G$. Its eigenvalues are called the Laplacian eigenvalues of $G$, which are denoted by $0=\mu_{1}(G)\leq\mu_{2}(G)\leq\cdots\leq \mu_{N}(G)$. $\mu_{2}(G)$ is called the {\em algebraic connectivity} of $G$.  As is well known, if $G$ is $d$-regular, then $\lambda_{1}(G)=d$ and $\mu_{2}(G)=d-\lambda_{2}(G)$. The lower bound of graph toughness in terms of eigenvalues are extensively studied (see \cite{Brouwer,Brouwer96,CLXHaemers,GZ,Gu,LK}).

Let $G$ be a connected and non-complete $d$-regular graph of order $N$. Set $b=-\lambda_{N}(G)$.
The well-known Hoffman ratio bound  \cite{Hoff} states that $\alpha(G)\leq\frac{bN}{d+b}$. If equality holds, one can observe that  $t(G)\leq\frac{d}{b}$ by considering a cut consisting of the complement of a maximum independent set of $G$.
Cioab\u{a} and Wong~\cite{CW} showed that $t(G)=\frac{d}{b}$ holds for several classes of regular graphs (attaining
equality in the Hoffman ratio bound), including complements of point graphs of generalized quadrangles, lattice graphs, Kneser graphs $\KG(n,2)$ (except the Petersen graph) and their complements.  Park et al.~\cite{Park} confirmed this for more Kneser graphs $\KG(n,k)$, where $k=3,4$ and $k\geq5$ for large $n$. Here, the vertices of the Kneser graph $\KG(n,k)$ are the $k$-subsets of $[n]=\{1,\ldots,n\}$, where $k\ge2$ and $n\ge2k+1$. Two vertices are adjacent if they are disjoint. It is shown that Kneser graphs (except the Petersen graph $\KG(5,2)$) have a Hamilton cycle in~\cite{MMN}.
Their work \cite{CW,Park} can be read as follows.

\begin{theorem}[\cite{CW}]\label{thm:k=2}
Let $n\ge6$ and $G=\KG(n,2)$. Then 
$t(G)=\frac{n-2}{2}$. Moreover, each tough set is precisely the complement of a maximum independent set in $G$.
\end{theorem}

For $G=\KG(5,2)$ (i.e., the Petersen graph), Cioab\u{a} and Wong~\cite{CW} also showed that $t(G)=\frac{4}{3}$. Moreover, a tough set is precisely a maximum independent set of size 4, leaving three disjoint edges.

\begin{theorem}[\cite{Park}]\label{thm:k=3,4}
Let $k=3,4,n\ge2k+1$ and $G=\KG(n,k)$. Then 
$t(G)=\frac{n-k}{k}$. Moreover, each tough set is precisely the complement of a maximum independent set in $G$.
\end{theorem}

 Park et al.~\cite{Park} also showed  $t(\KG(n,k))=\frac{n-k}{k}$ for $k\geq5$  and $n\geq\frac{2}{\ln2}k^{2}+(2-\frac{3}{\ln2})k+\frac{1}{\ln2}$. Moreover, each tough set is precisely the complement of a maximum independent set in $G$. They \cite{Park} proposed the following conjecture for other values of $n$.

\medskip

$\mathbf{Conjecture~ 1.}$ (\cite{Park})
Let $k\geq5,n\ge2k+1$ and $G=\KG(n,k)$. Then 
$t(G)=\frac{n-k}{k}$. Moreover, each tough set is precisely the complement of a maximum independent set in $G$.

\medskip

In this paper, we confirm Conjecture 1.

\begin{theorem}\label{thm:ordinary}
Conjecture 1 is correct.
\end{theorem}

Combining Theorem \ref{thm:ordinary} with the work \cite{CW,Park}, the toughness of Kneser graphs is completely determined, and the tough sets are also characterized. To end this section, we introduce an elementary result for Kneser graphs.
Let $G=\KG(n,k)$, where $k\ge2$ and $n\ge2k+1$.   For $x\in[n]$
the \emph{point-star} $\mathcal S_x=\{B\subseteq[n]:|B|=k,\ x\in B\}$
is an independent set of $G$. The Erd\H{o}s--Ko--Rado theorem and its
equality case~\cite{EKR,HM} identify the $\mathcal S_x$ as all the
largest intersecting families. Thus, $\alpha(G)=\binom{n-1}{k-1}$.

The rest of the paper is organized as follows. In Section 2, we include some elementary spectral results. In Section 3, we give some additional lemmas for our main result. The proof of Theorem \ref{thm:ordinary} is given in Section 4.

\section{Preliminaries}\label{sec:prelim}

In this section, we will include two elementary spectral results. Let $G$ be a graph. For $X\subseteq V(G)$, let $\partial X$ be the set of edges with exactly
one end vertex in $X$.

\begin{lemma}[\cite{Mo}]\label{lem:rayleigh}
Let $G$ be a connected graph of order $N$. For any $X\subseteq V(G)$, one has
$$|\partial X|\ge \frac{\mu_{2}(G)|X|(N-|X|)}{N}.$$ In particular, if
$h=\min\{|X|,N-|X|\}$, then $|\partial X|\ge\mu_{2}(G) h/2$.
\end{lemma}

The spectrum of the adjacency matrix of Kneser graphs were determined (see Chapter 9 of \cite{GR}). 

\begin{lemma}[\cite{GR}]\label{lem:spectrum}
For $k\geq2$ and $n\geq2k+1$, the eigenvalues of the adjacency matrix of $\KG(n,k)$ are
$(-1)^{j}\binom{n-k-j}{k-j}$ with multiplicities $\binom{n}{j}-\binom{n}{j-1}$ for $0\leq j\leq k$.
\end{lemma}

\section{Additional lemmas on components of a cut}

Let $G$ be a $d$-regular graph of order $N$. Set $b=-\lambda_{N}(G)$ and $\mu=d-\lambda_{2}(G)$. Let $S$ be a cut of $G$. For a component $C$ of $G-S$, the average internal degree of $C$ is defined as the average degree of $C$ as a graph itself. If $a$ denotes the average internal degree of $C$, define the weight of $C$ as
\[
       W(C)=\frac{|\partial C|}{b+a}.
\]
For $j\geq1$, let $t_j=j(d-j+1)/(b+j-1)$. Clearly, a clique of order
$j$ has average internal degree $j-1$, so its weight is
$t_j$ whenever it occurs as a component of $G-S$. The following two lemmas are helpful for our main result. We point out that the following Lemma \ref{lem:GZ} appeared in {\cite[eq.~(2)]{GZ}} in the Laplacian form. For the sake of completeness, we give the detailed proof here.

\begin{lemma}\label{lem:GZ}
Let $G$ be a $d$-regular graph of order $N$.  Set $b=-\lambda_{N}(G)$. Let $C_1,\ldots,C_c$ be the components of $G-S$ for some $S\subseteq V(G)$. Then $|S|\ge\sum_{1\leq i\leq c} W(C_{i})$.
\end{lemma}

\begin{proof} For $1\leq i\leq c$, set $m_i=|C_i|$, $a_i=2|E(C_i)|/m_i$ and
$e_i=|\partial C_i|$. It suffices to show
\[
 \qquad |S|\ \ge\ \sum_{i=1}^c \frac{e_i}{b+a_i}.
\]
Clearly, $e_i=m_i(d-a_i)$. Let $A$ denote the adjacency matrix of $G$.
 Take a vector $\mathbf{y}=(y_{v})_{v\in V(G)}$, such that equals $y_{v}=1$ for
$v\in S$ and  $y_{v}=-h_i$ for $v\in C_i$. Let $\mathbf{y}'$ denote the transpose of $\mathbf{y}$. Since $A+bI$ is positive
semidefinite, we obtain
\[
 0\le \mathbf{y}'(A+bI)\mathbf{y} =b|S|+2|E(S)|+
       \sum_{i=1}^c\bigl(m_i(b+a_i)h_i^2-2e_i h_i\bigr).
\]
Substituting
$h_i=e_i/[m_i(b+a_i)]$ and using $2|E(S)|=d|S|-\sum_i e_i$, we have 
\[
 (d+b)|S|\ge\sum_{1\leq i\leq c}\left(e_i+\frac{e_i^2}{m_i(b+a_i)}\right)
 =\sum_{1\leq i\leq c} e_i\,\frac{d+b}{b+a_i}.
\]
Hence, \[
 \qquad |S|\ \ge\ \sum_{i=1}^c \frac{e_i}{b+a_i}.
\]
\end{proof}

\begin{lemma}\label{lem:criterion}
Let $G$ be a connected noncomplete $d$-regular graph of order $N$.  Set $b=-\lambda_{N}(G)$ and $\mu=d-\lambda_{2}(G)$. Let $1\le j\le d$
be an integer and $t_j=j(d-j+1)/(b+j-1)$.  If
\begin{equation}\label{eq:level}
                    \mu(b-1)>2j(d+b),
\end{equation}
then $W(C)\ge t_j$ for every component $C$ with $|C|\geq j$ of $G-S$ (with $S$  a
cut). Equality holds precisely when $C$ is a
clique of order $j$.
\end{lemma}

\begin{proof}
Set $m=|C|$, $a=2|E(C)|/m$, $e=|\partial C|$, and
$T=d+1-j-t_j$. Clearly, $a\le m-1$ with equality
precisely when $C$ is complete. 
Hence, it suffices to show $e-t_j(b+a)\geq0$, with equality precisely when $C$ is a
clique of order $j$.
Noting $t_j(b+j-1)=j(d+1-j)$, we have
\begin{align}
 e-t_j(b+a)&=md-t_jb-(m+t_j)a \notag\\
 &\geq md-t_jb-(m+t_j)(m-1)\notag\\
 &=(m-j)(T-m).
 \label{eq:factor}
\end{align}
 For $m=j$
the right-hand side is zero and equality requires $a=j-1$ (implying that $C$ is a clique of order $j$). For
$j<m<T$ it is strictly positive. 

It remains to consider $m>j$ and $m\ge T$. Choose a vertex $y$ in
a different component from $C$ of $G-S$. Both $y$ and all its $d$ neighbours
lie outside $C$, so $N-m\ge d+1\ge T$. From \eqref{eq:level} we see $b>1$. 
Thus,
\[
     T=\frac{(d+1-j)(b-1)}{b+j-1}=\frac{(b-1)t_j}{j}>0.
\]
Using Lemma~\ref{lem:rayleigh} we have
\[
   e\ge\frac{\mu}{2}\min(m,N-m)\ge\frac{\mu T}{2}
      >t_j(d+b)\ge t_j(a+b),
\]
as desired.
\end{proof}

\section{Proof of Theorem \ref{thm:ordinary}}\label{sec:criterion}

Throughout this section, let $G=\KG(n,k)$, where $k\geq5$ and $n\geq2k+1$. Set $w=n-k$ and $r=\frac{w}{k}$. Let $d$ and $N$ denote the degree and the order of $G$, respectively. By Lemma \ref{lem:spectrum}, we have 
\[
    d=\binom wk,\qquad b=-\lambda_{N}(G)=\binom{w-1}{k-1},\qquad
    \mu=d-\lambda_2(G)=d-\binom{w-2}{k-2}.
\]

 Let $S$ be a cut of $G$. Recall that $W(C)=\frac{|\partial C|}{b+a}$ for a component $C$ of $S$, where $a$ is the average internal  degree of $C$. We first show the following lemma.

\begin{lemma}\label{important lem}
Let $C$ be a component of $G-S$ with $|C|\geq2$. Then
$W(C)\ge 2(d-1)/(b+1)$, with equality precisely when $|C|=2$.
\end{lemma}

\begin{proof} Set $R=2(d-1)/(b+1)$. It suffices to show $W(C)\ge R$, with equality precisely when $|C|=2$. Clearly, $W(C)=R$ when $|C|=2$. Without loss of generality, we can assume $|C|\geq3$ in the following. 
 We will prove this by two cases.

\medskip

$\mathbf{Case~ 1.}$ $n\ge2k+2$.

\medskip

Recall $w=n-k\geq k+2$. By Lemma~\ref{lem:spectrum}, we see $\lambda_2(G)/d= k(k-1)/(w(w-1))$. Noting  $d/(d+b)=r/(r+1)$, we obtain that
\[
 \frac{\mu(b-1)}{d+b}\ \ge\
 \frac{r}{r+1}(b-1)
 \left(1-\frac{k(k-1)}{w(w-1)}\right).
\]
Each factor on the right is positive and nondecreasing for $w\ge k+2$.
At $w=k+2$ their product is
\[
      \frac{2k+1}{2}\cdot\frac{(k-1)(k+2)}{(k+1)^2}=\frac{2k+1}{2}(1-1/(k+1)-2/(k+1)^2).
\]
This expression increases
with $k$, and at $k=5$ it is $77/18>4$. Thus, $\mu(b-1)>4(d+b)$. By Lemma \ref{lem:criterion}, we have $W(C)>R$ as $|C|\geq3$.

\medskip

$\mathbf{Case~ 2.}$ $n=2k+1$.

\medskip

In this case: $n=2k+1$, the parameters become $d=k+1$, $b=k$ and
$R=2k/(k+1)$. $G=\KG(2k+1,k)$ is triangle-free, since three pairwise
disjoint $k$-sets require $3k>n$ points. Also, the union of two
distinct $k$-sets contains at least $k+1$ points. Its complement
contains at most $k$ points. Thus, any two distinct vertices have at most one
common neighbour.

\medskip

Let $m=|C|\ge3$ and let $a$ be its average internal degree. We have two subcases.

\medskip

$\mathbf{Subcase~ 2.1.}$ $m\le2k-2$.

\medskip

Recall that the graph $C$ is triangle-free. Thus $\deg_C(v_{1})+\deg_C(v_{2})\leq m$ for any edge $v_{1}v_{2}$ of $C$.
 Summing over the edges
gives 
$$\sum_{v\in V(C)}\deg_C(v)^2=\sum_{vv'\in E(C)}(\deg_C(v)+\deg_C(v'))\le m|E(C)|.$$
 Using Cauchy's inequality, we have 
 $$4|E(C)|^{2}=\left(\sum_{v\in V(C)}\deg_C(v)\right)^{2}\leq m\sum_{v\in V(C)}\deg_C(v)^2.$$
It follows that $ma^2\le m^2a/2$, and thus 
$$a\le m/2.$$  
Recall $|\partial C|=m(d-a)$. Consequently, noting $a\leq \frac{m}{2}$, we have
\[
 |\partial C|-R(k+a)
  \ \ge\ m(k+1-m/2)-R(k+m/2)
  =\frac{(m-2)\{2k^2-(k+1)m\}}{2(k+1)} .
\]
The expression in braces is at least $2$ for $m\le2k-2$.
It follows that $W(C)>R$ as $m>2$.

\medskip

$\mathbf{Subcase~ 2.2.}$ $m\ge2k-1$.

\medskip

Recall $N=\binom{2k+1}{k}$. Set $h=\min(m,N-m)$. Note that $m\geq2k-1\geq k+2$ as $k\geq5$. Since a vertex of another component (except $C$) and its $k+1$
neighbours lie outside $C$, we have $N-m\geq k+2$. Hence, $$h\ge k+2.$$
 Since $b+a=k+a\le2k+1$, we see $W(C)\ge\frac{|\partial C|}{2k+1}$. It suffices to show $|\partial C|\geq4k$, since it will be that 
 $$W(C)\geq4k/(2k+1)>2k/(k+1)=R.$$
 
At last,  we aim to prove that $|\partial C|\ge4k$. 
For $k+2\le h\le2k-1$, let $H$ be the subgraph on the smaller
side of the boundary $\partial C$, and write $a_H$ for its average degree. Note $h=|H|$.
 Let $L$ be set of ordered vertices of form $(v_{1},v_{2},v_{3})$, where $v_{2}$ is adjacent to both $v_{1}$ and $v_{3}$. Recall that any two distinct vertices have at most one common
neighbour.
 We have
\[
 \sum_{v\in H}\deg_H(v)(\deg_H(v)-1)=|L|
                \ \le\ h(h-1).
\]
Using Cauchy's inequality, we can show 
$$ a_H(a_H-1)\ \le\ \frac{1}{h}\sum_{v\in H}\deg_H(v)(\deg_H(v)-1).$$
 Hence, $a_H(a_H-1)\leq h-1$.
 Since
$a_H(a_H-1)\ge5a_H-9$ by $(a_H-3)^2\ge0$, we obtain
$$a_H\le(h+8)/5.$$
 Thus
\[
       |\partial C|=h(k+1-a_H)\ \ge\ h(5k-3-h)/5 .
\]
Clearly, $h(5k-3-h)/5\geq(k+2)(5k-3-(k+2))/5=4k+1+(k-5)(4k+3)/5$ for $k+2\leq h\leq 2k-1$. Thus,
$|\partial C|\geq4k$. 

It remains to consider $h\ge2k$. Note that $x(1-x/N)$ increases for $0\leq x\leq N/2$. 
We claim that $N>8k^2$. To check this, it holds for $k=5$ as $\binom{11}{5}=462$.
When $k$ increases by one, $N$ increases by a factor
$2(2k+3)/(k+2)>3$, while $k^2$ increases by a factor
$(1+1/k)^2<2$. Hence, we proved $N>8k^2$.
By Lemmas~\ref{lem:rayleigh}
and~\ref{lem:spectrum}, noting $2k\leq h\leq\frac{N}{2}$, we have
\[
   |\partial C|\ \ge\ 2h(1-h/N)\geq4k(1-2k/N)
       \ =\ 4k-8k^2/N\ >\ 4k-1.
\]
Thus $|\partial C|\geq4k$, as desired.

\end{proof}

Now we are ready to prove Theorem ~\ref{thm:ordinary}.

\begin{proof}[Proof of Theorem~\ref{thm:ordinary}] Recall that $G=\KG(n,k)$, where $k\geq5$ and $n\geq2k+1$. Let $S$ be a cut of $G$ such that $t(G)=\frac{|S|}{c(G-S)}$.
Set $c=c(G-S)$. Let $\ell$ be the number of components $C$ with $|C|\geq2$ of $G-S$. Note that $W(C)=r$ for $|C|=1$.  By
Lemma \ref{important lem},  $W(C)\geq R$ for $|C|\geq2$. Clearly,
\[
 R-r=\frac{r(b-1)-2}{b+1}>0,
\]
since $b=\binom{w-1}{k-1}\ge k\ge5$ and $r>1$. Recall that  $\alpha(G)=\binom{n-1}{k-1}$. Thus, $t(G)\leq\frac{N-\alpha(G)}{\alpha(G)}=r$ by considering a cut consisting of the complement of a maximum independent set.
Using Lemma~\ref{lem:GZ},  we obtain
$$|S|\ge r(c-\ell)+R\ell.$$
 Thus $|S|\ge rc$, implying $t(G)\geq r$. Consequently, $t(G)=r$. This forces $\ell=0$, i.e., each component of $G-S$ is a singleton. Hence, $t(G)=\frac{N-c}{c}$. Note that $c\leq\alpha(G)$. This forces that $c=\alpha(G)$, since $t(G)=r=\frac{N-\alpha(G)}{\alpha(G)}$.  This completes the proof.

\end{proof}

\section*{Acknowledgements}

This work is supported by the National Natural Science Foundation of China (No.12301449, 12301448) and the Scientific Innovation Project for Young Scientists in Shandong Provincial Universities (No.2024KJG014).

\section*{Data availability}

No data were used for the research described in this article.

\section*{Competing Interests}
The authors declare that they have no known competing financial interests or
personal relationships that could have appeared to influence the work reported in
the paper.

\section*{Declaration on the Use of AI}
We used the OpenAI Codex (GPT-6)  for  proofreading, calculation checking, grammar checking, and language polishing.  The authors take full
responsibility for the correctness and originality of the paper.


\begin{thebibliography}{99}
\small\raggedright
\setlength{\itemsep}{2pt}




\bibitem{BauerHakimiSchmeichel1990}
D.~Bauer, S.~L. Hakimi, and E.~Schmeichel,
\emph{Recognizing tough graphs is NP-hard},
Discrete Appl. Math. \textbf{28} (1990), 191--195.

\bibitem{Brouwer} A. E. Brouwer, Toughness and spectrum of a graph,
\emph{Linear Algebra Appl.} \textbf{226--228} (1995), 267--271.

\bibitem{Brouwer96} A. E. Brouwer, Spectrum and connectivity of graphs,
\emph{CWI Quart.} \textbf{9} (1996), 37--40.



\bibitem{CLXHaemers} H. Chen, J. Li and S.-J. Xu,
A proof of Haemers's toughness conjecture,
\emph{European J. Combin.} \textbf{139} (2027), 104465.
\href{https://doi.org/10.1016/j.ejc.2026.104465}{doi:10.1016/j.ejc.2026.104465}.

\bibitem{Chvatal} V. Chv\'atal, Tough graphs and Hamiltonian circuits,
\emph{Discrete Math.} \textbf{5} (1973), 215--228.
\bibitem{CW} S. M. Cioab\u{a} and W. Wong, The spectrum and toughness of
regular graphs, \emph{Discrete Appl. Math.} \textbf{176} (2014), 43--52.

\bibitem{EKR} P. Erd\H{o}s, C. Ko and R. Rado, Intersection theorems for
systems of finite sets, \emph{Quart. J. Math. Oxford Ser. (2)}
\textbf{12} (1961), 313--320.


\bibitem{GR}
 C. Godsil and G. Royle, Algebraic Graph Theory, Graduate Texts in Mathematics, vol. 207, Springer-Verlag, New York, 2001.
 
\bibitem{GZ} G. Greaves and H. Zhu, A proof of Haemers' toughness conjecture,
preprint (2026). \href{https://arxiv.org/abs/2605.15738}{arXiv:2605.15738}.
\bibitem{Gu} X. Gu, A proof of Brouwer's toughness conjecture,
\emph{SIAM J. Discrete Math.} \textbf{35} (2021), 948--952.
\bibitem{HM} A. J. W. Hilton and E. C. Milner, Some intersection theorems
for systems of finite sets, \emph{Quart. J. Math. Oxford Ser. (2)}
\textbf{18} (1967), 369--384.

\bibitem{Hoff}
 A.J. Hoffman, On eigenvalues and colorings of graphs, in: Graph Theory and Its Applications, Proc. Advanced Sem., Math. Research Center, Univ. of
Wisconsin, Madison, Wis., 1969, Academic Press, New York, 1970, pp. 79-91.

\bibitem{LK} Y. Li and C. Y. Ku, Graphs attaining Haemers' toughness bound,
preprint (2026). \href{https://arxiv.org/abs/2611.18193}{arXiv:2611.18193}.


\bibitem{Mo}
 B. Mohar, Some applications of Laplace eigenvalues of graphs, in: Graph Symmetry, Springer, 1997, pp. 225-275.

\bibitem{MMN} A. Merino, T. M\"utze and Namrata, Kneser graphs are
Hamiltonian, \emph{Adv. Math.} \textbf{468} (2025), 110189.
\href{https://doi.org/10.1016/j.aim.2025.110189}{doi:10.1016/j.aim.2025.110189}.


\bibitem{Park} D. Park, A. Ostuni, N. Hayes, A. Banerjee, T. Wakhare,
W. Wong and S. Cioab\u{a}, The toughness of Kneser graphs,
\emph{Discrete Math.} \textbf{344} (2021), 112484.
\href{https://doi.org/10.1016/j.disc.2021.112484}{doi:10.1016/j.disc.2021.112484}.

\end{thebibliography}
\end{document}